\documentclass{ifacconf}

\usepackage{natbib}            
\usepackage{graphicx}          

\usepackage[update,prepend]{epstopdf}

\usepackage{amsmath} 
\usepackage{amssymb}  

\begin{document}

\begin{frontmatter}

\title{Energy-Efficient Building HVAC Control Using Hybrid System LBMPC
\thanksref{footnoteinfo}} 

\thanks[footnoteinfo]{This material is based upon work supported by the National Science Foundation under Grants CNS-0931843 (CPS-ActionWebs) and CNS-0932209 (CPS-LoCal). The views and conclusions contained in this document are those of the authors and should not be interpreted as representing the official policies, either expressed or implied, of the National Science Foundation.}

\author[First]{Anil Aswani} 
\author[First]{Neal Master} 
\author[First]{Jay Taneja}
\author[First]{Andrew Krioukov}
\author[First]{David Culler}
\author[First]{Claire Tomlin}

\address[First]{Electrical Engineering and Computer Sciences, University of California, Berkeley, CA 94720 USA (e-mail: \{aaswani,tomlin\}@eecs.berkeley.edu, neal.m.master@berkeley.edu, \{taneja,krioukov,culler\}@cs.berkeley.edu).}    
          

\begin{abstract}                          
Improving the energy-efficiency of heating, ventilation, and air-conditioning (HVAC) systems has the potential to realize large economic and societal benefits.  This paper concerns the system identification of a hybrid system model of a building-wide HVAC system and its subsequent control using a hybrid system formulation of learning-based model predictive control (LBMPC).  Here, the learning refers to model updates to the hybrid system model that incorporate the heating effects due to occupancy, solar effects, outside air temperature (OAT), and equipment, in addition to integrator dynamics inherently present in low-level control.  Though we make significant modeling simplifications, our corresponding controller that uses this model is able to experimentally achieve a large reduction in energy usage without any degradations in occupant comfort.  It is in this way that we justify the modeling simplifications that we have made.  We conclude by presenting results from experiments on our building HVAC testbed, which show an average of 1.5MWh of energy savings per day (p = 0.002) with a 95\% confidence interval of 1.0MWh to 2.1MWh of energy savings.
\end{abstract}

\end{frontmatter}

\section{Introduction}

Nearly 10\% of greenhouse gas emissions and 25\% of the electricity used in the United States is due to heating, ventilation, and air-conditioning (HVAC) systems in buildings \citep{BEDB2009, mcquade2009}.  This has driven research into better control methods (e.g., \citep{nghiem2011,kelman2011,frauke2012,liao2012,aswani2011_proc,ma2012}) that can help mitigate the negative externalities due to the large energy consumption of HVAC, while still ensuring the comfort of building occupants.  But the heterogeneity of HVAC equipment with respect to their physical modalities makes it difficult to develop a  control design methodology that scales to many types of equipment.

Even when new HVAC controllers are designed, experimentally comparing their efficiency and comfort in relation to existing controllers is difficult because of the large temporal variability in weather and occupancy conditions.  Identifying energy models of HVAC equipment can be difficult because some equipment is designed to operate most efficiently at certain temperatures or settings, which requires extensive measurement to characterize.  Moreover, not all buildings have equipment to directly measure the energy consumption of only HVAC equipment.

\subsection{Experiments with HVAC Controllers}

\citet{siroky2011} showed that a linear model predictive controller (MPC) could provide a 10-20\% reduction in energy usage of a ceiling radiant heating system, as compared to the default, manufacturer-provided controller.  The temperature difference between a water set point and the outside air temperature (OAT) was used to approximate the energy usage of each controller.


\citet{aswani2011_proc} designed a new controller for an air-conditioner on a single-room testbed, which achieved up to 30\% savings on warm days and up to 70\% savings on cool days.  Mathematical models of the temperature dynamics of different control schemes and their energy characteristics were constructed in order to allow comparisons between different controllers of experimentally measured HVAC energy usage to simulations over identical weather and occupancy conditions.

A 20\% improvement in performance of thermal storage for campus-wide building cooling was achieved by \citet{ma2012} using better control methods.  Direct energy measurements of the equipment along with a regression model of baseline performance were used to compare controllers.  This analysis approach requires direct measurements, because it can only statistically distinguish large differences.

\subsection{Overview}

This paper describes our design methodology for an energy-efficient controller of a building-wide HVAC system that is able to maintain comfortable occupant conditions.  We begin by describing our HVAC testbed, before explaining the modeling procedure (cf. \citep{aswani_2012hvac}) that was used to identify the thermal dynamics of the building and the HVAC system.  Next, we describe a hybrid system \citep{tomlin00} version of learning-based model predictive control (LBMPC) \citep{aswani2011_safe}.

LBMPC is a robust form of adaptive MPC.  Compared to linear parameter-varying MPC \citep{kothare2000,falcone2008}, LBMPC differs in that it provides robustness to model changes using tube MPC (e.g., \citep{chisci2001}).  Furthermore, the robust, adaptive MPC in \citep{Fukushima2007301,adetola2011} use an adaptive model with an uncertainty measure to ensure robustness, while LBMPC uses an adaptive model to improve performance and a nominal model with an uncertainty measure to provide robustness.


We conclude by experimentally comparing, on a building-wide HVAC system, our hybrid system LBMPC controller to the default controller.  This is done using the comparison methodology described in \citep{aswani2012_compare}, which uses nonparametric methods that can compute and compare quantitative metrics of energy usage and occupant comfort for different HVAC controllers.

\section{BRITE-S Testbed}

\begin{figure}
\begin{center}
\includegraphics[width=3.4in]{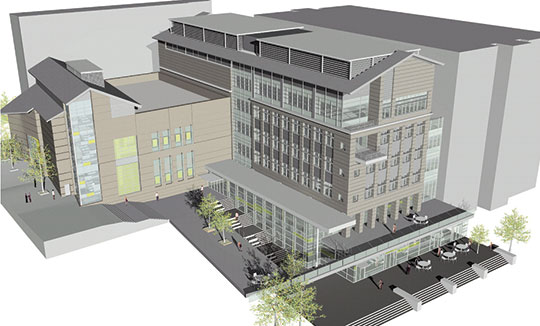}
\end{center}
\caption{\label{fig:sdh}Sutardja Dai Hall is 141,000-square-foot building, and is part of our BRITE-S testbed.}
\end{figure}

The Berkeley Retrofitted and Inexpensive HVAC Testbed for Energy Efficiency in Sutardja Dai Hall (BRITE-S) platform \citep{krioukov2011,aswani_2012hvac} is a building-wide HVAC system that maintains the indoor environment of a 141,000-square-foot building, shown in Fig. \ref{fig:sdh}, that is divided between a four-floor nanofabrication laboratory (NanoLab) and seven floors of general space (including office space, classrooms, and a coffee shop).  The building automation equipment can be measured and actuated through a BACnet protocol interface.

The HVAC system uses a 650-ton chiller to cool water.  Air-handler units (AHUs) with variable-frequency drive fans distribute air cooled by the water to variable air volume (VAV) boxes throughout the building.  Since the NanoLab must operate within tight tolerances, our control design can only modify the operation of the general space AHUs and VAV boxes, with no modification of chiller settings that are shared between the NanoLab and general space.

The default, manufacturer-provided controller in BRITE-S uses PID loops to actuate the VAV boxes and keeps a constant supply air temperature (SAT) within the AHUs. Conventional SAT reset control is not possible because several VAV boxes provide maximum air ﬂow rates throughout the day for nearly the entire range of SATs. These zones are often dominated by heating from computer equipment.

\section{Identifying Thermal Dynamics}

For the purpose of modeling and control, we will assume that data is sampled every 15 minutes; angle brackets (i.e., $\langle\cdot\rangle$) denote measurements sampled at this rate.  Let $T_s\langle k \rangle$ and $T_o \langle k \rangle$ be the SAT and OAT, respectively, at time $k$.  Similarly, $T_j\langle k \rangle$ is the temperature in the $j$-th zone of the building at time $k$, for $j = 1,\ldots,Z$ zones.  The VAV box in the $j$-th zone controls the zone temperature $T_j$ by modulating the amount of cool air sent to the zone $F_j\langle k \rangle$ and the amount that the air is reheated $R_j\langle k \rangle$.

In general, the thermal dynamics of each zone are
\begin{multline}
T_j\langle k+1 \rangle = f_j\Big(T_1\langle k \rangle, \ldots, T_Z\langle k \rangle, F_1\langle k \rangle,\ldots,F_Z\langle k \rangle, \\
R_1\langle k \rangle, \ldots, R_Z\langle k \rangle, T_o\langle k \rangle, T_s\langle k \rangle, O, X\Big),
\end{multline}
where $f_j(\cdot)$ is some unknown nonlinear function, $O$ are variables related to occupancy, and $X$ are variables related to other effects like the use of equipment, solar heating, etc.  Some simplifying assumptions are typically made by (a) considering the physical adjacency of different zones \citep{kelman2011,frauke2012,liao2012}, and (b) assuming that the effect of occupancy and other factors enters additively into the dynamics \citep{aswani2011_proc,aswani_2012hvac}.  After these assumptions, the model is
\begin{multline}
\label{eqn:nonlinear}
T_j\langle k+1 \rangle = f_j\Big(T_{n_1}\langle k \rangle, \ldots, T_{n_q}\langle k \rangle, F_j\langle k \rangle,R_j\langle k \rangle, \\
T_o\langle k \rangle, T_s\langle k \rangle\Big) + O + X,
\end{multline}
where $\{n_1,\ldots,n_q\}$ is the set of zones adjacent to $j$.

Additional assumptions allow further modeling simplifications.  We assume that the zone temperatures do not vary significantly throughout the day, since their temperature is in principle being controlled by the HVAC system.  Furthermore, we assume that the additive influence of the occupancy and other effects can be modeled by a single term $q_j\langle k \rangle$.  Even after making these assumptions, the model to be identified is nonlinear since the SAT $T_s$ affects the thermal dynamics in a bilinear form \citep{kelman2011,frauke2012}.


We take a hybrid system approach by forcing the SATs to belong to a finite set of values $\mathcal{M} = \{T_{s_1},\ldots,T_{s_p}\}$, where $p$ is the number of modes \citep{aswani_2012hvac}.  This allows us to consider multiple linearizations of (\ref{eqn:nonlinear}).  In our application to the BRITE-S testbed, we took $p = 3$ and $\mathcal{M} = \{52^\circ\text{F},58^\circ\text{F},62^\circ\text{F}\}$.  For the $m_i$-th mode (for $m_i \in \{1, \ldots, p\}$) of fixed SAT $T_{s_{m_i}}$, the model is given by
\begin{multline}
T_j\langle k+1 \rangle = a_{n_1}^{m_i}T_{n_1}\langle k \rangle + \ldots a_{n_q}^{m_i}T_{n_q}\langle k \rangle \\
+ b_j^{m_i}F_j\langle k \rangle + c_j^{m_i}R_j\langle k \rangle + d_j^{m_i} T_o\langle k \rangle + q_j\langle k \rangle,
\end{multline}
where the coefficients are unknown scalars and $q_j$ is an unknown function of time $k$.  The purpose of the system identification is to compute these unknown values.

\subsection{Experiments for System Identification}

In \citep{aswani_2012hvac}, we used one week of data with a fixed SAT to identify a model.  Identifying a hybrid system model where the SAT is able to change is more challenging, because identifying a model of equal fidelity would require three weeks (since we have three modes) of experimental data.  As a result, we modified our modeling approach.  We conducted experiments in which a small amount of data was gathered by cycling through all of the SATs in $\mathcal{M}$ so as to cover all of the modes of our hybrid system.  

Our experiment was as follows:  Starting at midnight, we set the SAT to $T_{s_1} = 52^\circ$F for two hours.  We next set the SAT to $T_{s_2} = 58^\circ$F for two hours.  After this, the SAT was changed to $T_{s_3} = 64^\circ$F for two hours.  During these six, consecutive hours, the other HVAC configuration was kept fixed.  The reason for picking a relatively quick horizon for all three experiments is that the heating load is roughly constant over a short time span.

\subsection{Initial Parameter Identification}


We used a small amount of training data to construct an approximate initial model that was used to do control.  To improve the controller performance, we have the option to re-identify the model using the semiparametric regression approach discussed in \citep{aswani_2012hvac}.

The approximate initial model was constructed as follows: We begin by making additional modeling simplifications.  Specifically, the exogenous heating load term $q_j\langle k \rangle$ was changed to also include the effects of OAT and adjacent zone temperatures.  For modeling purposes, we assume that the heating load term does not change significantly over a short time.  This was ensured by conducting the experiments for modeling within a quick time span.  We know \textit{a priori} that this last assumption is only approximate, but it serves to provide an initial model for which additional measurements can then be used to improve it.

Suppose we have (a) measurements for the $m_i$-th mode at times $k$ such that $L_{m_i} \leq k \leq U_{m_i}$; and (b) prior distributions for the coefficients $a_j^{m_i} \sim \mathcal{N}(\overline{a}_j^{m_i},\tilde{a}_j^{m_i}$), $b_j^{m_i} \sim \mathcal{N}(\overline{b}_j^{m_i},\tilde{b}_j^{m_i})$, $c_j^{m_i} \sim \mathcal{N}(\overline{c}_j^{m_i},\tilde{c}_j^{m_i})$, where the notation $\mathcal{N}(\mu,\Sigma)$ denotes a jointly Gaussian random variable with mean $\mu$ and covariance $\Sigma$.  Our initial model is given by
\begin{equation}
\label{eqn:dynamics}
T_j\langle k+1 \rangle = a_j^{m_i}T_j\langle k \rangle + b_j^{m_i}F_j\langle k \rangle + c_j^{m_i}R_j\langle k \rangle + q_j\langle k \rangle,
\end{equation}
and the coefficients can be identified by solving the following Bayesian, constrained least squares problem
\begin{align}
\min & \textstyle\sum_{m_i=1}^p\textstyle\sum_{k = L_{m_i}}^{U_{m_i}-1} (T_j\langle k+1 \rangle - a_j^{m_i}T_j\langle k \rangle \\
&\qquad - b_j^{m_i}F_j\langle k \rangle - c_j^{m_i}R_j\langle k \rangle - q_j\langle k \rangle)^2 \nonumber\\
&\qquad \quad + (a_j^{m_i} - \overline{a}_j^{m_i})^2/\tilde{a}_j^{m_i} \nonumber + (b_j^{m_i} - \overline{b}_j^{m_i})^2/\tilde{b}_j^{m_i} \\
&\qquad \qquad + (c_j^{m_i} - \overline{c}_j^{m_i})^2/\tilde{c}_j^{m_i} \nonumber\\
\text{s.t. } & a_j^r = a_j^s, \qquad\forall r,s \in \{1,\ldots,p\} \label{eqn:timeconstant}\\
& b_j^{r+1} < (T_{s_{r+1}}/T_{s_r})\cdot b_j^r, \qquad\forall r \in \{1,\ldots,p-1\} \label{eqn:cooling}\\
& c_j^r = c_j^s, \qquad\forall r,s \in \{1,\ldots,p\} \label{eqn:heating}\\
& q_j\langle k \rangle = q_j\langle q \rangle, \qquad\forall k,q \in \textstyle\bigcup_{r \in \{1,\ldots,p\}}[L_r,U_r]  \label{eqn:occupancy}
\end{align}

The constraints in the optimization problem reflect constraints between different hybrid modes of the HVAC system.  Constraint (\ref{eqn:timeconstant}) ensures that the time-constants of the thermal dynamics are constant for each mode of the hybrid system, while constraint (\ref{eqn:cooling}) encodes the fact that cooler temperatures will provide greater amounts of cooling in each zone.  Constraint (8) reflects that the re-heating capability of each VAV box is relatively constant across different SATs.  Lastly, constraint (\ref{eqn:occupancy}) represents the approximation that the occupancy is fixed over a short window of time.

\subsection{Modeling the VAV Box Control}

Each VAV box uses air flow $F_j$ and reheat amount $R_j$ to modulate the temperature of the $j$-th zone.  The VAV boxes use a PID controller (note that the use of a PID controller in zone control is typical for building HVAC systems \citep{honeywell}), and so we approximate this as a proportional controller.  This is not restrictive because the ``learning'' portion of our LBMPC controller can compensate for the unmodeled integrator portion of this control.  We will discuss this in the next section.

\begin{figure}
\begin{center}
\includegraphics{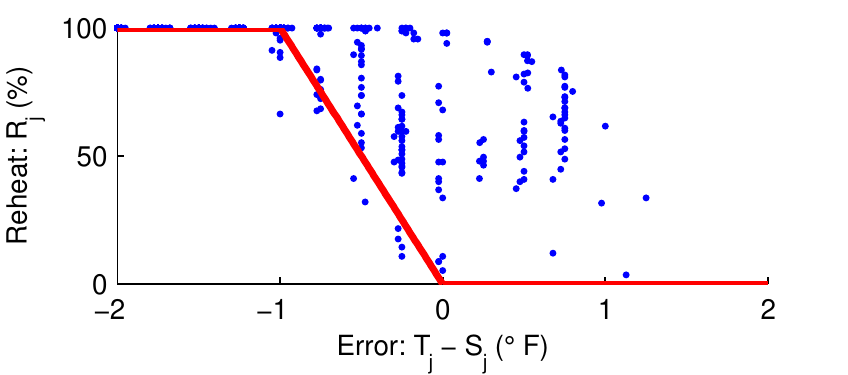}
\end{center}
\caption{\label{fig:reheat} A scatter plot of data taken from a single zone of temperature error and reheat amount in percentage is shown, along with the model (solid line) that we use for the VAV box controller.  The discrepancy is largely due to the integrator term that we leave unmodeled; our approach is to ``learn'' the value of the integrator term as we are doing control.}
\end{figure}

Let $e_j\langle k \rangle = T_j\langle k \rangle - S_j\langle k \rangle$ be the difference between the zone temperature and the temperature set point for the $j$-th zone ($S_j$).  Then, we use the following model for the control of reheat amount
\begin{equation}
\label{eqn:reheat}
R_j\langle k \rangle = \begin{cases} 100, & \text{if } e_j\langle k \rangle < -1^\circ\text{F} \\ -100\cdot e_j\langle k \rangle, & \text{if } -1^\circ\text{F} \leq e_j\langle k \rangle < 0^\circ\text{F} \\ 0, & \text{if } e_j\langle k \rangle \geq 0^\circ\text{F}
\end{cases}.
\end{equation}
We assumed that each VAV box has the same controller, and so this controller model was identified by fitting a piecewise linear model to the observed data points $(e_j\langle k \rangle, R_j\langle k \rangle)$ over all zones.  A comparison of the data points for a single zone and the fitted model is shown in Fig. \ref{fig:reheat}.

A similar model is used for the air flow amount
\begin{equation}
\label{eqn:flow}
F_j\langle k \rangle = \begin{cases} \alpha_j, & \text{if } e_j\langle k \rangle < 0^\circ\text{F} \\ (\omega_j - \alpha_j)\cdot e_j\langle k \rangle + \alpha_j, & \text{if } 0^\circ\text{F} \leq e_j\langle k \rangle < 1^\circ\text{F} \\ \omega_j, & \text{if } e_j\langle k \rangle \geq 1^\circ\text{F}
\end{cases},
\end{equation}
where $\alpha_j$ and $\omega_j$ are the minimum and maximum amount of air flow allowed in the $j$-th zone.  These values are configured by the building manager.

\section{The Hybrid System LBMPC Controller}

\label{sect:hl}
Several inputs can be used for control: The reheat $R_j$ and air flow $F_j$ in each zone can be explicitly actuated, or they can be implicitly actuated by varying the set point $S_j$.  Here, we only actuate the SAT; our methodology can be extended to also utilize the other inputs.


For the sake of argument, suppose that the sequence of SATs is fixed over a horizon of length $N$ and starting at time $k = 1$.  So if the mode sequence of SATs is $M = m_1,\ldots,m_N$, then the corresponding SATs are $T_s\langle 1 \rangle = T_{s_{m_1}},\ldots,T_s\langle N \rangle = T_{s_{m_N}}$.  Consequently, the temperature in each zone at time $k$ can be computed by just ``running'' or simulating the model forward in time; no optimization is needed to compute the values of $F_j$,$R_j$,$T_j$, etc.

As a result, minimizing our cost function subject to the thermal dynamics and system constraints requires optimization over \textit{only} the set of possible sequences.  In other words, our optimization problem is an integer program.  Fortunately, there is a constraint that greatly reduces the computational complexity: The SAT can only change once every hour; the reason for this constraint is that large, frequent changes to the SAT can damage the HVAC equipment.  So if the horizon is $N = 16$, which corresponds to four hours when sampling at 15 minute intervals, then there are 3621 different combinations.  

In order to further simplify the computational complexity, we make use of a heuristic that further reduces the computation.  We specify that the mode is fixed over the span of every hour.  For our setup --- where we have three modes and four possible mode changes --- this means that we have to compute the cost for $3^4 = 81$ different combinations.  This is a low number of combinations that can be computed under one second on a desktop computer because, as mentioned earlier, we do not need to optimize over other variables.  For future extensions where other variables are used to do control, we would only need to solve 81 convex optimization problems (specifically quadratic programs), which can be reasonably solved.

\subsection{Energy Modeling}

Before we present the optimization formulation of the controller, we provide some intuition into the form of the cost function we use.  A building HVAC system typically has several individual pieces of equipment that contribute to the overall energy usage.  Within BRITE-S, most of the energy consumption is due to three elements: the fans in the AHU, the chiller that cools water and indirectly cools air, and the reheating that occurs each zone's VAV box.

Even though parameterizations of the energy usage of the equipment are known \citep{kelman2011,frauke2012}, modeling these features is difficult because individual energy measurements are usually not available.  There is another subtle point regarding the energy models and their relationship to our cost function.  Let $E_1(\cdot),E_2(\cdot),E_3(\cdot)$ be the energy due to fans, chiller, and reheating.  One possible cost is $\textstyle\sum_j (T_j - S_j)^2 + \lambda E_1(\cdot) + \mu E_2(\cdot) + \gamma E_3(\cdot)$, where $\lambda,\mu,\gamma \in \mathbb{R}$ are constants.  There are two reasons for weighting energy usage based on type: First, the energy usage of the equipment is often related to the mechanical and physical stress on the equipment, and so differential weighting allows a finer level of regulation with respect to such considerations.  Second, the energy usage of equipment sometimes comes from different energy sources.  For example, heating is provided by steam while the fans are powered by electricity in BRITE-S.

Because the cost function weights each energy usage function, we really only need to know the energy model up to a constant, unknown scaling factor, because this gets subsumed into the scaling in the cost function.  This simplifies the energy modeling that we do for the purpose of control.  We use the following energy models: The fan energy usage is $E_1 \sim (\sum_j F_j)^3$ and the chiller energy is $E_2 \sim (T_o-T_s)\cdot\sum_j F_j$ \citep{kelman2011,frauke2012}, while the reheat energy is $E_3 \sim \sum_j R_j$.

\subsection{Optimization Formulation}

We can now present the optimization formulation of the hybrid system LBMPC controller.  The basic intuition behind LBMPC \citep{aswani2011_safe} is that two models of the system are kept.  The first is a nominal model that is used with respect to the constraints, and the second model is used in the cost function and is updated using data gathered during control.  This maintains robustness while improving performance through model updates.

Without loss of generality, we assume that the control action is being computed for time $k = 1$; also recall that we do control at a rate of every 15 minutes.  If $\overline{T}_j\langle 1 \rangle, \overline{R}_j\langle 1 \rangle, \overline{F}_j\langle 1 \rangle$ are the predictions of the linear model from time $k = 0$, then let $\hat{q}_j \langle i \rangle = T_j\langle 1 \rangle - \overline{T}_j \langle 1 \rangle$, $\hat{f}_j \langle i \rangle = F_j\langle 1 \rangle - \overline{F}_j \langle 1 \rangle$, $\hat{r}_j \langle i \rangle = R_j\langle 1 \rangle - \overline{R}_j \langle 1 \rangle$. The intuition is that these terms with hats represent corrections to the predictions of the MPC and provide the adaptation inherent in the controller \citep{aswani2011_proc}.  Specifically, $\hat{q}_j$ represents an estimate of the heating load due to occupants, weather, solar heating, and equipment; $\hat{f}_j$ represents the integrator term in the PID control of air flow in the $j$-th VAV box; and $\hat{r}_j$ represents the integrator term in the PID control of reheat amount in the $j$-th VAV box.

 The control action is given by the minimizer to
\begin{align}
\min_{m_1,\ldots,m_N} & \textstyle\sum_{i = 1}^N \bigg(\textstyle\sum_j (\tilde{T}_j\langle i+1 \rangle - S_j\langle i+1 \rangle)^2 \\
&\qquad \qquad + \lambda (\textstyle\sum_j \tilde{F}_j\langle i \rangle)^3 + \gamma \textstyle\sum_j \tilde{R}_j\langle i \rangle \nonumber\\
& \qquad \qquad \quad + \mu (T_o\langle i \rangle-T_s\langle i \rangle)\cdot\textstyle\sum_j \tilde{F}_j\langle i \rangle \bigg)\nonumber \\
\text{s.t.} & (\ref{eqn:dynamics}), (\ref{eqn:reheat}), (\ref{eqn:flow}) \nonumber\\
& \tilde{F}_j = F_j + \hat{f}_j;\ \tilde{R}_j = R_j + \hat{r}_j \nonumber \\
& \tilde{T}_j\langle i+1 \rangle = a_j^{m_i}\tilde{T}_j\langle i \rangle + b_j^{m_i}\tilde{F}_j\langle i \rangle + c_j^{m_i}\tilde{R}_j\langle i \rangle \nonumber \\
& \qquad \qquad \qquad + \hat{q}_j\langle i \rangle \nonumber \\
& T_s\langle i \rangle = T_{m_i} \nonumber \\
& T_s\langle 4q+r \rangle = T_s\langle 4q+s \rangle, \ \forall r,s \in \{1,\ldots,4\} \label{eqn:switching} \\
& \qquad \qquad \qquad \qquad \wedge q \in \{0,\ldots,\lceil N/4\rceil-1\} \nonumber \\
& T_j\langle i \rangle \in [66^\circ\text{F},78^\circ\text{F}], \ \forall j \in \{1,\ldots,Z\} \nonumber
\end{align}
Note that constraint (\ref{eqn:switching}) allows the SAT $T_s$ to switch value only once an hour, which reduces the computational complexity of the controller as discussed previously. 

A desktop computer took an average of under one second to solve this optimization problem for BRITE-S.  Furthermore, we used values of $\lambda = 6.7\text{e}4/(\sum_j \alpha_j)^3$, $\mu = 1.3\text{e$-$}3$, and $\gamma = 6.7$.  These values were picked using an iterative process in which (a) the control was computed but not used for actuation; (b) the control and predictions were analyzed for if SAT stayed at higher values, the reheat amount was low, and the air flow remained at moderate levels; (c) the coefficients $\lambda,\mu,\gamma$ were changed and the process starting at (a) was repeated until (b) was satisfied.

\section{Measuring Efficiency}

Separate energy measurements of the general space HVAC system are not available in BRITE-S.  Instead, we only have access to measurements of both the general space HVAC and the NanoLab HVAC.  We denote these energy measurements as $E[i]$ for $i = 1,\ldots,D$, where the index $i$ is over hourly intervals.  Furthermore, we have measurements of the OAT that correspond to the energy usage measurements: $T_o[i]$ for $i = 1,\ldots,D$.

The general model describing the relationship between energy usage, OAT, occupancy $O$, and other factors (e.g., solar effects, equipment, etc.) $X$ is $E = f(T_o, O, X)$, where $f(\cdot,\cdot,\cdot)$ is an unknown, nonlinear relationship.  But because occupancy and other factors are generally not directly measured, it is typically only possible to consider the energy usage averaged over occupancy and other factors $E_{o,x}(T_o) = \mathbb{E}\big[f(T_o, O, X) \big| T_o\big]$.  This can be estimated using nonparametric regression \citep{gyorfi2002}, and it is a curve that describes the relationship between the OAT and the average energy consumption.  Intuitively, if the data points $(T_o[i], E[i])$ for $i = 1,\ldots,D$ represent a scatter plot of energy usage versus OAT; then $\hat{E}_{o,x}(T_o)$ represents the smoothed version of the scatter plot.


The average amount of energy used in one hour is therefore $E_{g} = \mathbb{E}_{g}\big(E_{o,x}(T_o)\big) = \mathbb{E}_{g}\big[\mathbb{E}\big[f(T_o, O, X) \big| T_o\big]\big]$, where $g(T_o)$ is a probability distribution of OATs.  This notation allows us to define the average amount of energy used in one day as $E_{day} = \sum_{i = 1}^{24} E_{g_i}$, where $g_i(T_o)$ is the probability distribution of OAT during the $i$-th hour of the day.  For simplicity, we assume a uniform distribution for the OAT.


This allows us to compare the energy efficiency of two controllers.  We can compute the quantities defined above for each controller, and a nonparametric methodology \citep{boostrap2007} can be used to determine whether the differences are statistically significant.

\section{Measuring Occupant Comfort}
In order to define a quantification of comfort that is both tractable and will scale to many buildings, we focus on a measure that is only dependent on temperature.  We assume that the average temperature for the $j$-th zone is measured at hourly intervals $T_j[i]$ for $i = 1,\ldots,D$.



Let $(x)_+$ be the thresholding function, which is defined so that $(x)_+ = 0$ if $x < 0$ and $(x)_+ = x$ otherwise.  We define our quantification of comfort using soft thresholding as $C = 1/Z \cdot \textstyle\sum_{j=1}^Z \textstyle\int_0^1 (|T_j - S_j|-B_j)_+ dt$, where the integral with respect to $dt$ is over one hour of time, $S_j$ is the set point of the zone, and $B_j$ is the amount of temperature deviation for which the building has been configured.  The intuition is that this quantity increases whenever $T_j$ exceeds $S_j$ by more than $B_j$, and the amount of increase in this quantity is proportional to the amount and duration of temperature deviation.  The BRITE-S building is configured for $B_j \equiv 1^\circ$F for all zones and times.  


\section{Experimental Results}

The hybrid system LBMPC controller was used to control the SAT in BRITE-S for 8 days that spanned both weekdays and weekends, and this was compared to 22 days in which the default, manufacturer-provided controller was used.  Recall that our comparison methods implicitly account for variations in energy usage due to occupancy, and this is made explicit in the 95\% confidence interval for estimated values.  The building and HVAC configurations were kept identical for when both controllers were used, and this configuration has been in use for over one year.

For notational reasons, we use a superscript 1 to refer to the default controller, and superscript 2 denotes the hybrid system LBMPC controller.  Also, we define the quantity $\Delta\hat{E}^{2,1} = \hat{E}_{day}^2-\hat{E}_{day}^1$ to be the estimated difference between the average energy usage over a day of the LBMPC controller and that of the default controller.  The value $\Delta\hat{C}^{2,1} = \hat{C}_{day}^2-\hat{C}_{day}^1$ is an analogous quantity of the estimated difference in average comfort over a day.

\begin{figure}
\begin{center}
\includegraphics[clip = true, trim = 0in 0.10in 0in 0.21in]{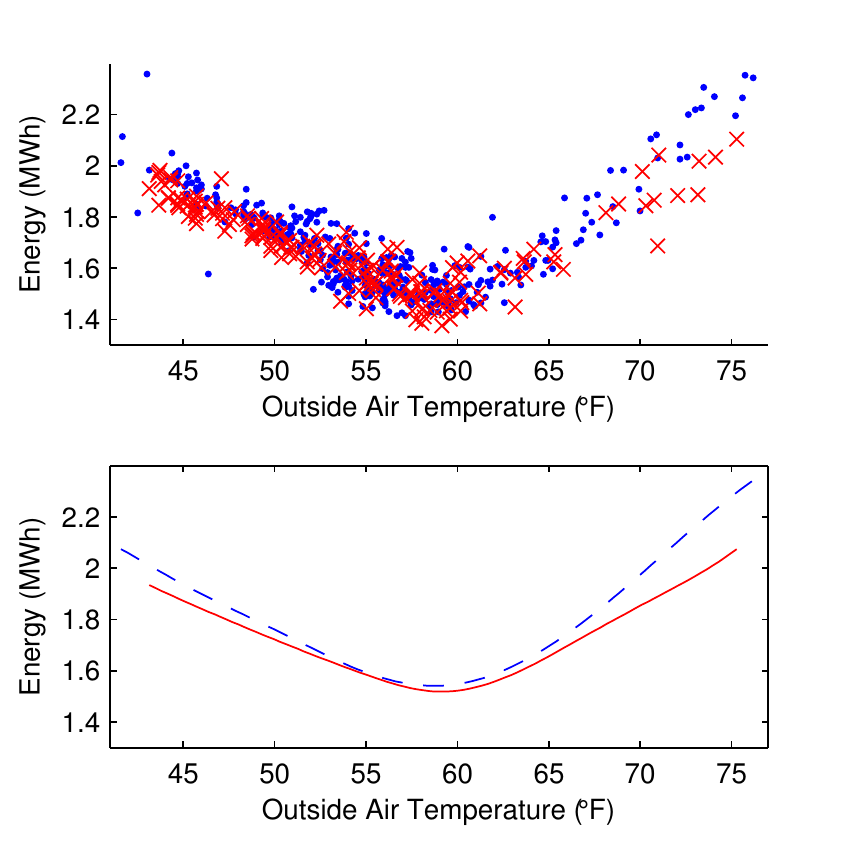}
\end{center}
\caption{\label{fig:lbmpc_energy}The cross marks and solid line (points and dashed line) denote the energy characteristics of the LBMPC (default) controller.}\end{figure}

The estimated energy characteristics $\hat{E}_{o,x}^1(T_0),\hat{E}_{o,x}^2(T_0)$ are shown in Fig. \ref{fig:lbmpc_energy}, and their differences are statistically significant $(p = 0.002)$.  The estimated difference in average energy usage over one day $\Delta \hat{E}^{2,1} = -1.53$MWh is statistically significant $(p=0.002)$.  And the 95\% confidence interval is $\Delta \hat{E}^{2,1} \in [-2.07,-1.02]$MWh.

\begin{figure}
\begin{center}
\includegraphics[clip = true, trim = 0in 0.10in 0in 0.21in]{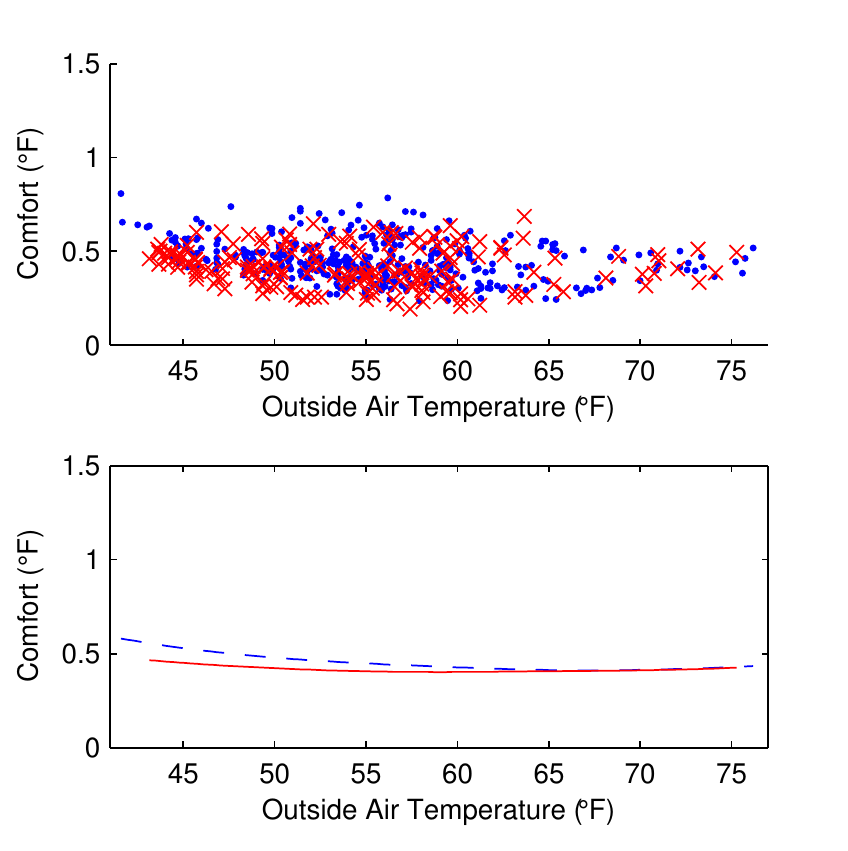}
\end{center}
\caption{\label{fig:lbmpc_comfort}The cross marks and solid line (points and dashed line) denote the comfort characteristics of the LBMPC (default) controller.}
\end{figure}

The estimated comfort characteristics $\hat{C}_{o,x}^1(T_0),\hat{C}_{o,x}^2(T_0)$ are shown in Fig. \ref{fig:lbmpc_comfort}, and their differences are not statistically significant $(p = 0.8)$.  Furthermore, the estimated difference in average comfort over one day $\Delta \hat{C}^{2,1} = -0.75^\circ$F is not statistically significant $(p=0.5)$, meaning that there is not enough evidence to exclude that $\Delta C^{2,1} = 0^\circ$F.

The LBMPC controller provides modest energy savings at most OATs, which sum up to significant savings over a day.  And because the difference in comfort characteristics is not statistically significant, this suggests that the LBMPC and default controllers provides comparable levels of comfort.

\section{Conclusion}

We have presented a hybrid model of building HVAC and described its control using hybrid system LBMPC.  Experiments show substantial savings, and future directions for further energy savings were discussed.  More broadly speaking, our experiments on BRITE-S, and previously on BRITE, show that the LBMPC methodology can provide significant energy savings for a wide variety of HVAC systems operating using different physical modalities.

\begin{ack}                               
The authors thank Domenico Caramagno and Stephen Dawson-Haggerty for their assistance. 
\end{ack}

\bibliography{ifacconf}           

\end{document}